\documentclass[10pt,a4paper]{article}
\usepackage{indentfirst}
\usepackage{amsmath,amssymb,amsfonts,amsthm,graphics}
\usepackage{epsfig}
\usepackage{makeidx}
\usepackage{mathrsfs}
\usepackage{color}
\usepackage{cite}
\usepackage{dsfont}
\usepackage{mathrsfs}
\usepackage{algorithmic}
\usepackage{algorithm}
\usepackage{graphicx}
\usepackage{geometry}
\usepackage{epstopdf}
\usepackage{setspace}
\usepackage{caption}
\usepackage{subfig}
\usepackage[section]{placeins}
\usepackage{float}
\newtheorem{observation}{Observation}

\begin{document}
\title{\Large\bf  A Note on  Improved bounds for the Oriented Radius of Mixed Multigraphs}

\author{Hengzhe Li$^{a}$, Zhiwei Ding$^{a}$, Jianbing Liu$^{b}$, Yanhong Gao$^{a}$, Shuli Zhao$^{a}$\\
\small $^{a}$College of Mathematics and Information Science,\\
\small Henan Normal University, Xinxiang 453007, P.R. China\\
\small $^b$ Department of Mathematics, \\
\small University of Hartford, West Hartford 06117, USA\\
\small Email: lihengzhe@htu.edu.cn, dzhiweidzw@163.com,\\ \small jianliu@hartford.edu, gaoyanhong@htu.edu.cn, zhaoshuli0210@126.com}
\date{}
\maketitle
\begin{abstract}

For a positive integer $r$, let $f(r)$ denote the smallest number such that any 2-edge connected mixed graph with radius $r$ has an oriented radius of at most $f(r)$. Recently, Babu, Benson, and Rajendraprasad significantly improved the upper bound of $f(r)$ by establishing that $f(r) \leq 1.5r^2 + r + 1$, see [Improved bounds for the oriented radius of mixed multigraphs, J. Graph Theory, 103 (2023), 674-689]. Additionally, they demonstrated that if each edge of a graph $G$ is contained within a cycle of length at most $\eta$, then the oriented radius of $G$ is at most $1.5r\eta$. The authors' results were derived through Observation 1, which served as the foundation for the development of Algorithm ORIENTOUT and Algorithm ORIENTIN. By integrating these algorithms, they obtained the improved bounds. However, an error has been identified in Observation 1, necessitating revisions to Algorithm ORIENTOUT and Algorithm ORIENTIN. In this note, we address the error and propose the necessary modifications to both algorithms, thereby ensuring the correctness of the conclusions.
\\[2mm]
\end{abstract}

\section{Introduction}
In this note, all graphs are considered to be finite mixed graphs. We refer to \cite{Babu} for any undefined notation and terminology in the following.

A {\it mixed multigraph} ({\it mixed graph} for brevity) $G$ is defined as an ordered pair $G=(V,E)$, where $V$ represents a set of vertices, and $E$ represents a multiset of unordered and ordered pairs of vertices. These pairs are referred to as undirected and directed edges, respectively.

A {\it walk} in a mixed graph $G$ from a vertex $x$  to a vertex $y$ is a sequence of vertices $v_1, v_2,\ldots, v_k$ such that $x=v_1$, $y=v_k$, and for every $1\le i\le k-1$, there exists either an undirected
edge $v_iv_{i+1}$ or an edge directed $v_iv_{i+1}$ from $v_i$ to $v_{i+1}$ in $G$.
A {\it trail} is defined as a walk that contains no repeated edges. A {\it path} is a trail that contains no repeated vertices, except possibly $v_1=v_k$.
For a path or cycle $C$, and for $x,y\in V(C)$, we use $C(x,y)$ to denote the subpath of $C$ from $x$ to $y$.

A {\it bridge} in a connected mixed multigraph $G$  is an undirected edge of $G$  whose removal disconnects the underlying undirected multigraph of $G$. A connected mixed graph that contains no bridges is {\it $2$-edge connected}. An {\it orientation} of a mixed graph $G$ is an assignment of exactly one direction to each undirected edge of $G$. A mixed graph $G$ is said to be {\it strongly orientable} if it can be oriented such that the resulting directed graph is strongly connected.

Given an undirected, directed, or mixed graph $G$, and a vertex $u\in V(G)$,
let $N_G(u)$ denote the set comprising of all its in‐neighbours, out‐neighbours, and undirected neighbours, and let $N_G[u]=N_G(u)\cup \{u\}$.
The {\it length} of a path in $G$ is defined as the number of edges (directed and undirected) contained in the path. For two distinct vertices $u$ and $v$, the {\it distance}
$d_G(u,v)$ in $G$ is the length of the shortest path from $u$ to $v$ in $G$. A subset $D\subseteq V(G)$ is an {\it $r$-step dominating set} of $G$ if there exists a vertex $u$ such that $d_G(u,v)\le r$ for each vertex $v\in V(G)\setminus D$.

The {\it out-eccentricity}  $e_{out}(u)$ of $u$ is defined as $\max\{d_G(u,v)|v\in V(G)\}$, and the {\it in-eccentricity}  $e_{in}(u)$ of $u$ is defined as $\max\{d_G(v,u)|v\in V(G)\}$. The {\it eccentricity}  $e(u)$ of $u$ is $\max\{e_{out}(u), e_{in}(u)\}$. The {\it radius} $rad(G)$ (respectively, {\it diameter} $diam(G)$) of $G$ is $\min\{e(u)| u\in V(G)\}$ (respectively, $\max\{e(u)| u\in V(G)\}$). A vertex $u$ of $G$ is referred to as a {\it center vertex} if $e(u)=rad(G)$.

The {\it oriented radius}
(respectively, {\it oriented diameter}) of $G$ is defined as the minimum radius (respectively, diameter) of an orientation of $G$. For each $r\in\mathds{N}$, let $f(r)$  (respectively, $\bar{f}(r)$) represent the smallest value such that each $2$-edge connected mixed graph (respectively, undirected graph) with radius $r$ possesses oriented radius at most $f(r)$ (respectively, $\bar{f}(r)$). For each $d\in\mathds{N}$, let $g(d)$ (respectively, $\bar{g}(d)$) represent the minimum value such that any $2$-edge connected mixed
graph (respectively, undirected graph) with diameter $d$ has oriented diameter at most $g(d)$ (respectively, $\bar{g}(d)$).

In 1939, Robbins in \cite{Robbins} solved the One-Way Street Problem and proved that a graph $G$ admits a strongly connected orientation if and only if $G$ is bridgeless, that is, $G$ does not have any cut-edge. Naturally, one hopes that the oriented diameter of a bridgeless graph is as small as possible. In 1978, Chv\'atal and Thomassen \cite{Chvatal} showed that $\bar{f}(r)=r^2+r$ and $\frac{1}{2}d^2+d\le\bar{g}(d)\le 2d^2+2d$. One can see that $8\le\bar{g}(3)\le 24$.
In 2010, Kwok, Liu and West \cite{Kwok} narrowed the gap between the upper and the lower bounds by showing that $9\le \bar{g}(3)\le 11$.
In 2022, Wang and Chen determined that $\bar{g}(3)=9$.

For mixed graph,  Chung, Garey, and Tarjan \cite{Chung} in 1985 showed that $r^2+r\le f(r)\le 4r^2+4r$. Recently, Babu, Benson, and Rajendraprasad made a big improvement on the upper bound by showing that $f(r)\le 1.5r^2+r+1$. Moreover, they showed that if each edge of $G$ lies in a cycle of length at most $\eta$, then the oriented radius of $G$ is at most $1.5r\eta$.
In their paper, authors presented Observation~$1$.
Based on this observation, they developed Algorithm ORIENTOUT and Algorithm ORIENTIN, and by combining these two algorithms, they derived the aforementioned results. However, there is an error in Observation $1$, which necessitates adjustments to Algorithm ORIENTOUT and Algorithm ORIENTIN.

In this note, we rectify the error and make the necessary adjustments to Algorithm ORIENTOUT and Algorithm ORIENTIN to ensure the validity of the conclusions.
\section{Main Result}
In this section, we fix the error in Observation~$1$, and adjust {\bf Algorithm} ORIENTOUT and {\bf Algorithm} ORIENTIN in \cite{Chung}.

\noindent\rule{\linewidth}{0.4pt}
{\bf Algorithm} ORIENTOUT

\noindent\rule{\linewidth}{0.4pt}
Input: A strongly connected bridgeless mixed multigraph $G$ and a vertex $u \in V(G)$ with eccentricity at most $r$.

\noindent Output: An orientation $\overrightarrow{H}$ of a subgraph $H$ of $G$ such that $N[u]\subseteq V(\overrightarrow{H})$
and for every vertex $v$ in $ \overrightarrow{H}$, $ d_{\overrightarrow{H}}(u,v) \leqslant 2r$ and  $d_{\overrightarrow{H}}(v,u) \leqslant 4r-1$.

\noindent\rule{\linewidth}{0.4pt}

We create $\overrightarrow{H}$  in four stages, starting from Stage 0.

Stage $0$.  Let $v$ be a vertex having multiple edges incident with $u$. If possible, these edges are oriented in such a way that all the $uv$ edges are part of a directed $2$-cycle. Notice that this does
not increase any pairwise distances in $G$. We can also remove multiple oriented edges in the same direction between any pair of vertices without affecting any distance. Hence the only multiedges left in $G$ are those which form a directed $2$-cycle. We denote the resulting graph as ${G_0}$.

Let $X$ denote the set of all vertices with at least one edge incident with the vertex $u$. $X$ is partitioned into $X_{in}$, $X_{out}$, and $ X_{un}$.
A vertex $v \in X$ is said to be in $X_{in}$ if it has at least one directed edge towards $u$.
A vertex $v \in (X\setminus X_{in})$ is said to be in $X_{out}$ if it has at least one directed edge from $u$.
 Finally, $X_{un}=X\setminus(X_{in}\cup X_{out})$. Notice that a vertex $v \in X_{un}$ has exactly one undirected edge incident with $u$.
 We initialize $X_{conf}=\emptyset$ (we will later identify this as the set of conflicted vertices in $X$).
 For each $v \in X$, let $l(v)$ denote the length of the shortest cycle containing an edge between $u$ and $v$,
 $l(v) \leqslant  (2r + 1)$ by Lemma~3 in \cite{Babu}.
 Let $ s=\sum_{v \in X} l(v)$ .

Stage 1. Orient some of the undirected edges of $G_{0}$ incident with $u$ in
this stage to obtain a mixed graph $G_1$ as follows.

Repeat Steps (i)-(iii) for all the vertices in $X_{un}$ :

(i) An edge $uv_i$ is oriented from $v_i$ to $u$ if the parameter $s$ remains the same even after such an orientation. The vertex $v_i$ is added to $X_{in}$ in this case.

(ii) Otherwise, the edge $uv_i$ is oriented from $u$ to $v_i$ if the parameter $s$ remains the same even after such an orientation. The vertex $v_i$ is added to $X_{out}$ in this case.

(iii) Otherwise, we leave the edges $uv_i$ unoriented. Such an edge $uv_i$ is called a {\it conflicted edge} and the vertex $v_i$ is called a {\it conflicted vertex}. The vertex $v_i$ is added to $X_{conf}$ in this case.

\begin{observation}\label{Bug}{\rm (Babu, Benson, and Rajendraprasad \cite{Babu})} If an edge $uv_i$  is conflicted then there exists an edge $\overrightarrow{v_ju}$, at the time of
processing $v_i$, where $v_j\in  X_{in}$ and $j<i$ such that the edge $uv_i$ is a part of every shortest
cycle containing the edge $\overrightarrow{v_ju}$. Otherwise the parameter $s$ would have remained the same
even if the edge $uv_i$ gets oriented from $u$ to $v_i$. Hence, for every vertex $v\in X_{conf}$, $uv$ is an
undirected edge and there exists a $w\in X_{in}$ such that every shortest path from $u$ to $w$ starts with the edge $uv$.
\end{observation}

It is unfortunate that there is a bug in Observation \ref{Bug}. See the counterexample in Figure~$1$. It is easy to check that $X_{in}=\{v_1,v_2\}$, $X_{out}=\{v_3,v_4\}$, and $X_{conf}=\{v_5\}$. For $v_5\in X_{conf}$, by Observation \ref{Bug}, there is a vertex $v_j\in X_{in}$ such that the edge $uv_5$ is a part of every shortest cycle containing the edge $\overrightarrow{v_ju}$. But, for $j=1,2$,  every shortest cycle containing the edge $\overrightarrow{v_ju}$ does not contain the edge $uv_5$, a contradiction.

\begin{center}
\scalebox{0.85}[0.85]{\includegraphics{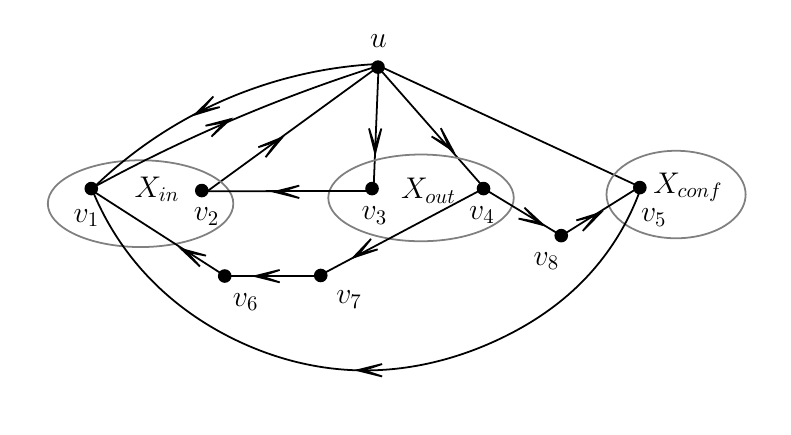}}

\small
Figure 1. A counterexample for Observation~$1$.
\end{center}

We can fix Observation \ref{Bug} as follows.

\begin{observation}\label{Fix} {\rm (Corrected version)}
If an edge $uv_i$ is conflicted,
then either there exists a vertex $v_j \in  X_{in}$ such that the edge $uv_i$ is a part of every shortest cycle containing the edge $\overrightarrow{v_ju}$,
or there exists a vertex $v_k\in X_{out}$ such that the edge $uv_i$ is a part of every shortest cycle containing the edge $\overrightarrow{uv_k}$.
If not, that is, there are no such vertices in $X_{in}\cup X_{out}$, then there exists an orientation of the edge $uv_i$ such that the parameter $s$ remains the same. So, we can add $v_i$ into either $X_{in}$ or $X_{out}$, which contradicts to $v_i\in X_{conf}$. Hence, for every vertex $v\in X_{conf}$ and undirected edge $uv$, either there exists a $w \in X_{in}$ such that every shortest path from $u$ to $w$ starts with the edge $uv$, or there exists a $z \in X_{out}$ such that every shortest path from $z$ to $u$ ends with the edge $vu$.
\end{observation}

Stage $2$. Let $T_{out}$ be a minimal tree formed by a breadth‐first search in $G_1$, starting from $u$ and including a shortest path from $u$ to $v$ for each $v\in  X_{in}$.
Let $T_{in}$ be a minimal tree formed by a breadth‐first search in $G_1$ starting from $u$ and including a shortest path from $v$ to $u$ for each $v\in  X_{out}$. By Observation~\ref{Fix}, $X_{conf}\subseteq V(T_{out})\cup V(T_{in})$. We orient $T_{out}$ such that it is an out-tree rooted at $u$, and let $P_i$ be the only path in $T_{out}$ from $u$ to $v_i$ for each $v_i\in X_{in}$. The resulting graph is denoted by $G_2$.

For each $v_j\in X_{out}$, let $P_j$ denote the only path in $T_{in}$ from $v_j$ to $u$. Define $x_j$ and $y_j$ as the first vertex and last vertices, respectively, in $V(P_j)\cap V(T_{out})$ when traversing from $v_j$ to $u$. 
If $V(P_j)\cap V(T_{out})=\emptyset$, then let $x_j=y_j=u$.
Let $C_j$ be the cycle formed by adding a directed edge $uv_j$ to $P_j$. Orient $C_j(v_j,x_j)$ to create a directed path from $v_j$ to $x_j$, and orient $C_j(y_j,u)$ to create a directed path from $y_j$ to $u$. Let $G_3$ be the result graph.

Stage 3. Every vertex in $V(T_{out})\setminus \{u\}$ is part of a directed cycle in $G_2$ that contains $u$ and has length at most $2r+1$.
Hence, for any vertex $x\in V(T_{out})\setminus \{u\}$,
we have $d_{G_2}(u,x) \le 2r $ and  $d_{G_2}(x,u)\le 2r$.

For each $y\in \bigcup_{v_j\in X_{out}}(V(C_j(u,x_j))\cup V(C_j(y_j,u)))$, there exists a vertex $v_j\in X_{out}$, such that $y\in V(P_j)$.
If $V(P_j)\bigcap (V(T_{out})\setminus \{u\})=\emptyset$, then $y$ belongs to a directed cycle in $G_3$ that contains $u$ and has length at most $2r+1$, and so $d_{G_3}(u,y) \leqslant 2r$ and  $d_{G_3}(y,u) \leqslant 2r$. Otherwise, that is, $V(P_j)\bigcap (V(T_{out})\setminus \{u\})\neq \emptyset$, then $d_{G_3}(y,(V(T_{out})\setminus \{u\}))\le 2r-1$, which implies that $d_{G_3}(y,u)\le 4r-1$. 

Considering $d_{G_3}(u,y)$. If $y\in V(C_j(u,x_j))$, then $d_{G_3}(u,y)\le 2r$. Suppose that $y\in V(C_j(y_j,u))$, and $y_j\in V(P_i)\cup V(P_j)$, where $P_i$ is the shortest path from $u$ to $v_i$ in $T_{out}$.
Then $P_i(u,y_j)\cup C_j(y_j,y)$ cotains a path from $u$ to $y$ with length at most $|E(P_i(u,y_j))|+|E(C_j(y_j,y))|\le |E(C_j(u,y_j))|+|E(C_j(y_j,y))|\le 2r$, since $|E(P_i(u,y_j))|\le |E(C_j(u,y_j))|$ by the definition $T_{out}$ and $P_i$. Let $H=G_3[V(T_{out})\cup (\bigcup_{v_j\in X_{out}}(V(C_j(u,x_j))\cup V(C_j(y_j,u))))]$. For each undirected edge in $H$, we can orient it arbitrarily, and let $\overrightarrow{H}$ be the result graph. It is easy to check that $d_{\overrightarrow{H}}(u,v)\le 2r$ and $d_{\overrightarrow{H}}(v,u)\le 4r-1$ for each $v\in V(\overrightarrow{H})$.

For {\bf Algorithm} ORIENTIN, we also need to fix it similarly. After the above fix, the same result holds, see \cite{Babu} for details.

\section{Acknowledgments}
This work was supported by the Postdoctoral Research Grant in Henan Province under Grant No. HN2022148 and
Natural Science Foundation of Henan Province under Grant No. 242300420645.

\end{document}